\documentclass[12pt]{amsart}
\usepackage{amsmath, amsthm, amscd, amsfonts}

\newtheorem{theorem}{Theorem}[section]
\newtheorem{lemma}[theorem]{Lemma}
\newtheorem{proposition}[theorem]{Proposition}
\newtheorem{corollary}[theorem]{Corollary}
\theoremstyle{definition}

\newtheorem{example}[theorem]{Example}
\theoremstyle{remark}

\numberwithin{equation}{section}
\begin{document}
\title[Vanishing of the first $(\sigma,\tau)$-cohomology group]{Vanishing of the first $(\sigma,\tau)$-cohomology group
of triangular Banach algebras}
\author[M. Khosravi, M.S. Moslehian and A.N. Motlagh]{M. Khosravi$^1$, M. S. Moslehian$^2$ and A. N. Motlagh$^3$}
\address{$^1$Maryam Khosravi: Department of Mathematics, Teacher Training
University, Tehran, Iran; \newline Banach Mathematical Research
Group (BMRG), Mashhad, Iran.} \email{khosravi$_-$m@saba.tmu.ac.ir}
\address{$^2$Mohammad Sal Moslehian: Department of Mathematics, Ferdowsi University, P. O. Box 1159, Mashhad 91775, Iran;\newline Centre of Excellence in Analysis on Algebraic
Structures (CEAAS), Ferdowsi Univ., Iran.}
\email{moslehian@ferdowsi.um.ac.ir}
\address{$^3$Abolfazl Niazi Motlagh: Department of Mathematics, Ferdowsi University, P. O. Box 1159, Mashhad 91775, Iran.}
\email{ab$_{_{-}}$ni40@stu-mail.um.ac.ir and niazimotlagh@gmail.com}
\subjclass[2000]{Primary 46H25; Secondary 46L57, 16E40.}
\keywords{Banach algebra, triangular Banach algebra, Banach
bimodule, $(\sigma,\tau)$-derivation, first
$(\sigma,\tau)$-cohomology group, $(\sigma,\tau)$-amenability,
$(\sigma,\tau)$-weak amenability.}
\begin{abstract}
In this paper, we define the first topological
$(\sigma,\tau)$-cohomology group and examine vanishing of the first
$(\sigma,\tau)$-cohomology groups of certain triangular Banach
algebras. We apply our results to study the $(\sigma,\tau)$-weak
amenability and $(\sigma,\tau)$-amenability of triangular Banach
algebras.
\end{abstract}
\maketitle

\section{Introduction and preliminaries}

Suppose that $\mathcal A$ and $\mathcal B$ are two unital
algebras with units $1_{\mathcal A}$ and $1_{\mathcal B}$,
respectively. Recall that a vector space $\mathcal M$ is a unital
$\mathcal A-\mathcal B$-bimodule whenever it is both a left
$\mathcal A$-module and a right $\mathcal B$-module satisfying
$$a(mb)=(am)b, 1_{\mathcal A}m=m1_{\mathcal B}=m \quad (a, b \in\mathcal A, m\in\mathcal M).$$
Then ${\rm Tri}(\mathcal A,\mathcal M,\mathcal B)=\left [
\begin{array}{cc}\mathcal A&\mathcal M\\0&\mathcal B \end{array}\right ]=\{ \left [
\begin{array}{cc}a&m\\0&b \end{array}\right ]; a\in\mathcal A, m\in\mathcal M,
b\in\mathcal B \}$ equipped with the usual $2\times 2$
matrix-like addition and matrix-like multiplication is an algebra.

An algebra ${\mathcal T}$ is called a triangular algebra if there
exist algebras $\mathcal A$ and $\mathcal B$ and nonzero $\mathcal
A-\mathcal B$-bimodule $\mathcal M$ such that ${\mathcal T}$ is
(algebraically) isomorphic to ${\rm Tri}(\mathcal A,\mathcal
M,\mathcal B)$. For example, the algebra ${\mathcal T}_n$ of
$n\times n$ upper triangular matrices over the complex field
${\mathbb C}$, may be viewed as a triangular algebra when $n>1$. In
fact, if $n>k$, we have ${\mathcal T}_n={\rm Tri}({\mathcal
T}_{n-k},M_{n-k,k}({\mathbb C}),{\mathcal T}_k)$ in which
$M_{n-k,k}({\mathbb C})$ is the space of $(n-k)\times k$ complex
matrices, cf. \cite{CHE}.

Let ${\mathcal T}$ be a triangular algebra. If $1=\left [
\begin{array}{cc}u&p\\0&v\end{array}\right ]$, and $\left [
\begin{array}{cc}a&0\\0&b \end{array}\right ]$ is denoted by
$a\oplus b$, then it can be easily verified that $e=u\oplus 0$ is an
idempotent such that $(1-e){\mathcal T}e=0$ but $e{\mathcal
T}(1-e)\neq 0$. Conversely, if there exists an idempotent $e\in
{\mathcal T}$ such that $(1-e){\mathcal T}e=0$ but $e{\mathcal
T}(1-e)\neq 0$. Then the mapping $x\mapsto \left [
\begin{array}{cc}exe&ex(1-e)\\0&(1-e)x(1-e)\end{array}\right ]$
is an isomorphism between ${\mathcal T}$ and ${\rm Tri}(e{\mathcal
T}e,e{\mathcal T}(1-e),(1-e){\mathcal T}(1-e))$; cf. \cite{CHE}.

By a triangular Banach algebra we mean a Banach algebra A which
is also a triangular algebra. Many algebras such as upper
triangular Banach algebras \cite{F-M1}, nest algebras \cite{DAV},
semi-nest algebras \cite{D-F-M}, and joins \cite{G-S} are
triangular algebras.

Following \cite{CHE}, consider a triangular Banach algebra
${\mathcal T}$ with an idempotent $e$ satisfying $e{\mathcal
T}(1-e)=0$ and $(1-e){\mathcal T}e\neq 0$. Put $\mathcal
A=e{\mathcal T}e, \mathcal B=(1-e){\mathcal T}(1-e)$ and $\mathcal
M=(1-e){\mathcal T}e$. Then $\mathcal A$ and $\mathcal B$ are
closed subalgebras of ${\mathcal T}$, $\mathcal M$ is a Banach
$\mathcal A-\mathcal B$-bimodule, and ${\mathcal T}={\rm
Tri}(\mathcal A,\mathcal M,\mathcal B)$. Conversely, given Banach
algebras $(\mathcal A,\|.\|_{\mathcal A})$ and $(\mathcal
B,\|.\|_{\mathcal B})$ and an $\mathcal A-\mathcal B$-bimodule
$\mathcal M$, then the triangular algebra ${\mathcal T}={\rm
Tri}(\mathcal A,\mathcal M,\mathcal B)$ is a Banach algebra with
respect to the norm given by $\|\left [
\begin{array}{cc}a&m\\0&b
\end{array}\right ]\|_{\mathcal T}=\|a\|_{\mathcal A}+\|m\|_{\mathcal M}+\|b\|_{\mathcal B}$.
It is not hard to show that each norm $\|.\|$ making ${\rm
Tri}(\mathcal A,\mathcal M,\mathcal B)$ into a triangular Banach
algebra is equivalent to $\|.\|_{\mathcal T}$, if the natural
restrictions of $\|.\|$ to $\mathcal A,\mathcal B$ and $\mathcal
M$ are equivalent to the given norms on $\mathcal A,\mathcal B$
and $\mathcal M$, respectively. See also \cite{M-S-Y, MOS1}

The concept of topological cohomology arose from the problems
concerning extensions by H. Kamowitz \cite{KAM}, derivations by R.
V. Kadison and J. R. Ringrose \cite{K-R1, K-R2} and amenability by
B.E. Johnson \cite{JOH} and has been extensively developed by A. Ya.
Helemskii and his school \cite{HEL1}. The reader is referred to
\cite{HEL1, RUN} for undefined notation and terminology.

Let $\mathcal A$ be a Banach algebra and $\sigma,\tau$ be
continuous homomorphisms on ${\mathcal  A}$. Suppose that
$\mathcal E$ is a Banach $\mathcal A$-bimodule. A linear mapping
$d:\mathcal A \to \mathcal E$ is called a
$(\sigma,\tau)$-derivation if
\begin{eqnarray*}
d(ab)=d(a)\sigma(b)+\tau(a)d(b)\quad (a, b \in\mathcal A).
\end{eqnarray*}
For example (i) Every ordinary derivation of an algebra $\mathcal
A$ into an $\mathcal A$-bimodule is an $id_{\mathcal
A}$-derivation, where $id_{\mathcal A}$ is the identity mapping on
the algebra $ \mathcal A$. (ii) Every point derivation $d:\mathcal A \to {\mathbb
C}$ at the character $\theta$ on $\mathcal A$ is a $\theta$-derivation.

A linear mapping $d:\mathcal A\longrightarrow \mathcal E$ is called
$(\sigma,\tau)$-inner derivation if there exists $x \in\mathcal E$
such that $d(a)=\tau(a)x-x\sigma(a) \quad (a \in\mathcal A)$. See
also \cite{M-M1, M-M2, MOS2, MOS3} and references therein.

We denote the set of continuous $(\sigma,\tau)$-derivations from
$\mathcal A$ into $\mathcal E$ by $Z^1_{(\sigma,\tau)}(\mathcal A,\mathcal E)$ and the set of
inner $(\sigma,\tau)$-derivations by $B^1_{(\sigma,\tau)}(\mathcal A,\mathcal E)$. we
define the space $H^1_{(\sigma,\tau)}(\mathcal A,\mathcal E)$ as the quotient space
$Z^1_{(\sigma,\tau)}(\mathcal A,\mathcal E)/B^1_{(\sigma,\tau)}(\mathcal A,\mathcal E)$.
The space $H^1_{(\sigma,\tau)}(\mathcal A,\mathcal E)$ is called the first
$(\sigma-\tau)$-cohomology group of $\mathcal A$ with coefficients in $\mathcal E$.

From now on, $\mathcal A$ and $\mathcal B$ denote unital Banach
algebras with units $1_{\mathcal A}$ and $1_{\mathcal B}$,
$\mathcal M$ denotes a unital Banach $\mathcal A-\mathcal
B$-bimodule and $\mathcal T={\rm Tri}(\mathcal A,\mathcal
M,\mathcal B)$ is the triangular matrix algebra. In addition,
$\mathcal X$ is a unital Banach ${\mathcal T}$-bimodule,
$\mathcal X_{\mathcal A\mathcal A}=1_{\mathcal A}\mathcal
X1_{\mathcal A}, \mathcal X_{\mathcal B\mathcal B}=1_{\mathcal
B}\mathcal X1_{\mathcal B}, \mathcal X_{\mathcal A\mathcal B}=
1_{\mathcal A}\mathcal X1_{\mathcal B}$ and $\mathcal X_{\mathcal
B\mathcal A}=1_{\mathcal B}\mathcal X1_{\mathcal A}$. For
instance, with $\mathcal X={\mathcal T}$ we have $\mathcal
X_{\mathcal A\mathcal A}=\mathcal A, \mathcal X_{\mathcal
B\mathcal B}=\mathcal B, \mathcal X_{\mathcal A \mathcal
B}=\mathcal M$ and $\mathcal X_{\mathcal B\mathcal A}=0$.

In this paper, we examine vanishing of the first
$(\sigma,\tau)$-cohomology groups of certain triangular Banach
algebras. We apply our results to investigate the
$(\sigma,\tau)$-weak amenability and $(\sigma,\tau)$-amenability of
triangular Banach algebras.

\section{Vanishing of the first $(\sigma,\tau)$-cohomology group}

In this section, using some ideas of \cite{F-M2}, we investigate
the relation between the first $(\sigma,\tau)$-cohomology  of
${\mathcal T}$ with coefficients in $\mathcal X$ and those of $\mathcal A$ and $\mathcal B$
with coefficients in $\mathcal X_{\mathcal A\mathcal A}$ and $\mathcal X_{\mathcal B\mathcal B}$,
respectively,  whenever $\mathcal X_{\mathcal A\mathcal B}=0$ in a direct method.

We start our work by investigating the structure of bounded
$(\sigma,\tau)$-derivations from a triangular Banach algebra into
bimodules.

Let $\sigma$ and $\tau$ be two homomorphisms on ${\mathcal T}$
with the following properties.
\begin{eqnarray}
\tau(1\oplus0)=1\oplus0,\, \tau(0\oplus1)=0\oplus1\label{equ1};\\
\sigma(1\oplus0)=1\oplus0, \,
\sigma(0\oplus1)=0\oplus1\label{equ2}.
\end{eqnarray}
The above relation implies easily that $\sigma(\mathcal
A)\subseteq\mathcal A$ and $\sigma(\mathcal B)\subseteq\mathcal B$
if we identify $a\in\mathcal A$ with $\left [
\begin{array}{cc}a&0\\0&0 \end{array}\right ]$  and $b\in\mathcal B$ with $\left [
\begin{array}{cc}0&0\\0&b \end{array}\right ]$. So with no ambiguity, we can consider $\sigma$ and $\tau$
as homomorphisms on $\mathcal A$ or $\mathcal B$, when it is
necessary. Now let $m\in \mathcal M$. If $\sigma\Big(\left [
\begin{array}{cc}0&m\\0&0 \end{array}\right ]\Big)=\left [
\begin{array}{cc}a'&m'\\0&b' \end{array}\right ]$, then
\begin{eqnarray*}
\left [\begin{array}{cc}a'&m'\\0&b' \end{array}\right
]=\sigma\Big(\left [
\begin{array}{cc}0&m\\0&0 \end{array}\right ]\Big)&=&\sigma\Big(\left [
\begin{array}{cc}1&0\\0&0 \end{array}\right ]\left [
\begin{array}{cc}0&m\\0&0 \end{array}\right ]\Big)\\
&=&\sigma\Big(\left [
\begin{array}{cc}1&0\\0&0 \end{array}\right ]\Big)\sigma\Big(\left [
\begin{array}{cc}0&m\\0&0 \end{array}\right ]\Big)\\
&=&\left [
\begin{array}{cc}1&0\\0&0 \end{array}\right ]\left [
\begin{array}{cc}a'&m'\\0&b' \end{array}\right ]\\
&=&\left [
\begin{array}{cc}a'&m'\\0&0 \end{array}\right ]
\end{eqnarray*}
and
\begin{eqnarray*}
\left [\begin{array}{cc}a'&m'\\0&b' \end{array}\right
]=\sigma\Big(\left [
\begin{array}{cc}0&m\\0&0 \end{array}\right ]\Big)&=&\sigma\Big(\left [
\begin{array}{cc}0&m\\0&0 \end{array}\right ]\left [
\begin{array}{cc}0&0\\0&1 \end{array}\right ]\Big)\\
&=&\sigma\Big(\left [
\begin{array}{cc}0&m\\0&0 \end{array}\right ]\Big)\sigma\Big(\left [
\begin{array}{cc}0&0\\0&1 \end{array}\right ]\Big)\\
&=&\left [
\begin{array}{cc}a'&m'\\0&b' \end{array}\right ]\left [
\begin{array}{cc}0&0\\0&1 \end{array}\right ]\\
&=&\left [
\begin{array}{cc}0&m'\\0&b' \end{array}\right ].
\end{eqnarray*}
Hence $\sigma\Big(\left [
\begin{array}{cc}0&\mathcal M\\0&0 \end{array}\right ]\Big)\subseteq\left [
\begin{array}{cc}0&\mathcal M\\0&0 \end{array}\right ]$. Thus one
can define $\sigma_{\mathcal M}:\mathcal M \to \mathcal M$ by
$m\mapsto m'$. To simplify the notation we denote $\sigma_{\mathcal
M}$ by $\sigma$. Thus $\sigma\Big(\left [
\begin{array}{cc}a&m\\0&b \end{array}\right ]\Big)$ can be written
as $\left [
\begin{array}{cc}\sigma(a)&\sigma(m)\\0&\sigma(b) \end{array}\right
]$.

\noindent If $\sigma_{\mathcal A}:\mathcal A \to \mathcal A$ and $
\sigma_{\mathcal B}:\mathcal B \to \mathcal B$ are homomorphisms,
then $\sigma_{\mathcal A}\oplus\sigma_{\mathcal B}:\mathcal
A\oplus\mathcal B \to \mathcal A\oplus\mathcal B$ defined by
$(\sigma_{\mathcal A}\oplus\sigma_{\mathcal
B})(a,b)=(\sigma_{\mathcal A}(a),\sigma_{\mathcal B}(b))$ is a
homomorphism. Conversely  every homomorphism on $\mathcal
A\oplus\mathcal B$ is of the form $\sigma_{\mathcal
A}\oplus\sigma_{\mathcal B}$ for some homomorphisms
$\sigma_{\mathcal A}$ and $\sigma_{\mathcal B}$ on $\mathcal A$ and
$\mathcal B$, respectively.

Applying our notation, let $\delta:{\mathcal T}\to\mathcal X$ be a bounded
$(\sigma,\tau)$-derivation. Then
$\delta_{\mathcal A}:\mathcal A \to 1_{\mathcal A}\mathcal X1_{\mathcal A}$ defined
by   $$\delta_{\mathcal A}(a)=1_{\mathcal A}\delta(\left [ \begin{array}{cc}a&0\\0&0
\end{array}\right ])1_{\mathcal A},$$ and
$\delta_{\mathcal B}:\mathcal B \to 1_{\mathcal B}\mathcal X1_{\mathcal B}$
defined by $$\delta_{\mathcal B}(b)=1_{\mathcal B} \delta(\left [
\begin{array}{cc}0&0\\0&b \end{array}\right ])1_{\mathcal B}$$ are bounded
$(\sigma,\tau)$-derivations.

Moreover, the mapping $\theta:\mathcal M \to 1_{\mathcal A}\mathcal X1_{\mathcal B}$ given by
$$\theta(m)=1_{\mathcal A}\delta(\left [ \begin{array}{cc}0&m\\0&0
\end{array}\right ])1_{\mathcal B}$$ satisfies
\begin{eqnarray}\label{2.1}
\theta(am)&=&1_{\mathcal A}\delta(\left [\begin{array}{cc}0&am\\0&0
\end{array}\right ])1_{\mathcal B}\nonumber\\
&=&1_{\mathcal A}\delta(\left [
\begin{array}{cc}a&0\\0&0 \end{array}\right ] \left [
\begin{array}{cc}0&m\\0&0 \end{array}\right ])1_{\mathcal B}\nonumber\\
&=&1_{\mathcal A}\tau(a)\delta(\left [
\begin{array}{cc}0&m\\0&0 \end{array}\right ])1_{\mathcal B}
+1_{\mathcal A}\delta(\left [ \begin{array}{cc}a&0\\0&0 \end{array}\right
])\sigma(\left [ \begin{array}{cc}0&m\\0&0 \end{array}\right
])1_{\mathcal B}\nonumber\\
&=&\tau(a)1_{\mathcal A}\delta(\left [
\begin{array}{cc}0&m\\0&0 \end{array}\right ])1_{\mathcal B}
+1_{\mathcal A}\delta(\left [ \begin{array}{cc}a&0\\0&0 \end{array}\right ])1_{\mathcal A}
\sigma(\left [ \begin{array}{cc}0&m\\0&0 \end{array}\right ])\nonumber\\
&=&\tau(a)\theta(m)+\delta_{\mathcal A}(a)\sigma(m)
\end{eqnarray}
and
\begin{eqnarray}\label{2.2}
\theta(mb)=\theta(m)\sigma(b)+\tau(m)\delta_{\mathcal B}(b).
\end{eqnarray}
Conversely, if $\delta_1$ and $\delta_2$ are bounded
$(\sigma,\tau)$-derivations of $ \mathcal A$ and $\mathcal B$ into
$\mathcal X_{\mathcal A\mathcal A}$ and $\mathcal X_{\mathcal
B\mathcal B}$, respectively, and $\theta:\mathcal M \to \mathcal
X_{\mathcal A\mathcal B}$ is any continuous linear mapping
satisfies (\ref{2.1}) and (\ref{2.2}), then the mapping $D(\left [
\begin{array}{cc}a&m\\0&b \end{array}\right
])=\delta_1(a)+\delta_2(b)+\theta(m)$ defines a bounded
$(\sigma,\tau)$-derivation of ${\mathcal T}$ into $X$, since
\begin{eqnarray*}
&&\tau(\left [ \begin{array}{cc}a&m\\0&b \end{array}\right ])
D(\left[
\begin{array}{cc}a'&m'\\0&b' \end{array}\right ])+D(\left [
\begin{array}{cc}a&m\\0&b \end{array}\right ])\sigma(\left [
\begin{array}{cc}a'&m'\\0&b' \end{array}\right ])\\
&=&\tau( \left [
\begin{array}{cc}a&m\\0&b \end{array}\right])
(\delta_1(a')+\delta_2(b')+\theta(m'))
+(\delta_1(a)+\delta_2(b)+\theta(m))\sigma(\left [
\begin{array}{cc}a'&m'\\0&b' \end{array}\right ])\\
&=&\tau(\left [
\begin{array}{cc}a&m\\0&b \end{array}\right ])\tau(1_{\mathcal A})\delta_1(a')
+\delta_1(a) \sigma(1_{\mathcal A})\sigma(\left [
\begin{array}{cc}a'&m'\\0&b'
\end{array}\right ])\\
&&+\tau(\left [
\begin{array}{cc}a&m\\0&b \end{array}\right ])\tau(1_{\mathcal B})\delta_2(b')
+\delta_2(b)\sigma(1_{\mathcal B})\sigma(\left [
\begin{array}{cc}a'&m'\\0&b'
\end{array}\right ])\\
&&+\tau(\left [
\begin{array}{cc}a&m\\0&b \end{array}\right ])\tau(1_{\mathcal A})
\theta(m')+\theta(m)\sigma(1_{\mathcal B})\sigma(\left [
\begin{array}{cc}a'&m'\\0&b'
\end{array}\right])\\
&=&\tau(a)\delta_1(a')+\delta_1(a)\sigma(a')+\delta_1(a)\sigma(m')+\tau(b)
\delta_2(b')\\
&&+\delta_2(b)\sigma(b')+\tau(m)\delta_2(b')+\tau(a)
\theta(m')+\theta(m)\sigma(b')\\
&=& \delta_1(aa')+\delta_2(bb')+\theta(am')+\theta(mb')\\
&=&D(\left [ \begin{array}{cc}aa'&am'+mb'\\0&bb' \end{array}\right ])\\
&=&D(\left [ \begin{array}{cc}a&m\\0&b \end{array}\right ]\left [
\begin{array}{cc}a'&m'\\0&b' \end{array}\right ]).
\end{eqnarray*}
If $\mathcal X_{\mathcal A\mathcal B}=0$, then we may assume that the linear mapping $\theta$
defined above is zero. Notice that, in this case,
$\delta_{\mathcal A}(a)\sigma(m)=\tau(m)\delta_{\mathcal B}(b)=0$ for every
$a\in\mathcal A, b\in\mathcal B, m\in\mathcal M$.

We are now ready to provide our main theorem.

\begin{theorem}\label{main}
Let $\mathcal X_{\mathcal A\mathcal B}=1_{\mathcal A}\mathcal X1_{\mathcal B}=0$.
Then $$H^1_{(\sigma,\tau)}({\mathcal T},\mathcal X)=H^1_{(\sigma,\tau)}(\mathcal A,\mathcal X_{\mathcal A\mathcal A})\oplus
H^1_{(\sigma,\tau)}(\mathcal B,\mathcal X_{\mathcal B\mathcal B}).$$
\end{theorem}
\begin{proof} Suppose that $ \mathcal X_{\mathcal A\mathcal B}=0$ and consider
the linear mapping $$\rho: Z^1_{(\sigma,\tau)}({\mathcal
T},\mathcal X) \to H^1_{(\sigma,\tau)}(\mathcal A,\mathcal
X_{\mathcal A\mathcal A})\oplus H^1_{(\sigma,\tau)}(\mathcal
B,\mathcal X_{\mathcal B\mathcal B})$$ defined by
$$\delta \mapsto (\delta_{\mathcal A}+N^1_{(\sigma,\tau)}(\mathcal A,\mathcal X_{\mathcal A\mathcal A}),
\delta_{\mathcal B}+N^1_{(\sigma, \tau)}(\mathcal B,\mathcal
X_{\mathcal B\mathcal B})).$$ If $\delta_1\in Z^1_{(\sigma,
\tau)}(\mathcal A,\mathcal X_{\mathcal A\mathcal A})$ and
$\delta_2\in Z^1_{(\sigma,\tau)}(\mathcal B,\mathcal X_{\mathcal
B\mathcal B})$, then $D(\left [
\begin{array}{cc}a&m\\0&b
\end{array}\right ])=\delta_1(a)+\delta_2(b)$ is a $(\sigma, \tau)$-derivation
from ${\mathcal T}$ into $\mathcal X$ and
\begin{eqnarray*}
\rho(D)&=&(D_{ A}+N^1_{(\sigma, \tau)}(\mathcal A,\mathcal
X_{\mathcal A\mathcal A}), D_{\mathcal B}+
N^1_{(\sigma, \tau)}(\mathcal B,\mathcal X_{\mathcal B\mathcal B}))\\
&=&(\delta_1+N^1_{(\sigma, \tau)}(\mathcal A,\mathcal X_{\mathcal
A\mathcal A}), \delta_2+ N^1_{(\sigma, \tau)}(\mathcal B,\mathcal
X_{\mathcal B\mathcal B}).
\end{eqnarray*}
The last equation is deduced from the fact that
$$D_{\mathcal A}(a)=1_{\mathcal A}(\delta_1(a)+\delta_2(0))1_{\mathcal A}=\delta_1(a),$$
and
$$\delta_{\mathcal B}(b)=1_{\mathcal B}(\delta_1(0)+\delta_2(b))1_{\mathcal B}=\delta_2(b).$$
Thus $\rho$ is surjective.

If $\delta \in {\rm ~ker}\rho$, then $\delta_{\mathcal A}\in
N^1_{(\sigma, \tau)}(\mathcal A,\mathcal X_{\mathcal A \mathcal
A})$ and $\delta_{\mathcal B} \in N^1_{(\sigma, \tau)}(\mathcal
B,\mathcal X_{\mathcal B\mathcal B})$. Then $\delta_{\mathcal
A}(a)=\tau(a)x-x\sigma(a)$ for some $x\in\mathcal X_{\mathcal
A\mathcal A}$ and $\delta_{\mathcal B}(b)=\tau(b)y-y\sigma(b)$ for
some $y\in\mathcal X_{\mathcal B\mathcal B}$. Then
\begin{eqnarray*}
D(\left [ \begin{array}{cc}a&m\\0&b \end{array}\right
])&=&\delta_{\mathcal A}(a)+
\delta_{\mathcal B}(b)\\
&=&(\tau(a)x-x\sigma(a))+(\tau(b)y-y\sigma(b))\\
&=&(\tau(a)+\tau(m)+\tau(b))(x+y)-(x+y)(\sigma(a)+\sigma(m)+\sigma(b))\\
&=&\tau(\left [ \begin{array}{cc}a&m\\0&b \end{array}\right ])(x+y)-(x+y)\sigma(\left [ \begin{array}{cc}a&m\\0&b \end{array}\right ]).
\end{eqnarray*}
Thus $D\in N^1_{(\sigma,\tau)}({\mathcal T},\mathcal X)$.\\
It is straightforward to show that
\begin{eqnarray*}
\delta(\left [ \begin{array}{cc}a&0\\0&0 \end{array}\right
])&=&1_{\mathcal A}\delta(\left [ \begin{array}{cc}a&0\\0&0
\end{array}\right ])1_{\mathcal A}+1_{\mathcal B}\delta(\left [ \begin{array}{cc}a&0\\0&0
\end{array}\right ])1_{\mathcal A}+1_{\mathcal B}\delta(\left
[ \begin{array}{cc}a&0\\0&0 \end{array}\right ])1_{\mathcal B}\\
&=& 1_{\mathcal A}\delta(\left [ \begin{array}{cc}a&0\\0&0
\end{array}\right ])1_{\mathcal A}+1_{\mathcal B}\delta(\left [
\begin{array}{cc}1_{\mathcal A}&0\\0&0
\end{array}\right ])1_{\mathcal A}\sigma(a).
\end{eqnarray*}
Similarly,
$$\delta(\left [
\begin{array}{cc}0&0\\0&b\end{array}\right ])=1_{\mathcal B}\delta(\left [
\begin{array}{cc}0&0\\0&b\end{array}\right
 ])1_{\mathcal B}-\tau(b)
 1_{\mathcal B}\delta(\left [ \begin{array}{cc}1_{\mathcal A}&0\\0&0\end{array}\right
  ])1_{\mathcal A},$$
and also
\begin{eqnarray*}
\delta(\left [
\begin{array}{cc}0&m\\0&0\end{array}\right ])&=&\delta(\left [
 \begin{array}{cc}1_{\mathcal A}&0\\0&0\end{array}\right ]\left[\begin{array}{cc}
 0&m\\0&0\end{array}\right ])\\
 &=& 1_{\mathcal B}\delta(\left [ \begin{array}{cc}1_{\mathcal A}&0\\0&0\end{array}\right
 ])1_{\mathcal A}
\sigma(\left[\begin{array}{cc}0&m\\0&b\end{array}\right])-\tau\left
[\begin{array}{cc}a&m\\0&0\end{array}\right])1_{\mathcal
B}\delta(\left [
\begin{array}{cc}1_{\mathcal A}&0\\0&0\end{array}\right ])1_{\mathcal A}.
\end{eqnarray*}
These follow that
\[\aligned
(\delta-D)(\left [ \begin{array}{cc}a&m\\0&b\end{array}\right
])&=\delta(\left [ \begin{array}{cc}a&m\\0&b\end{array}\right
])-1_{\mathcal A}\delta(\left [
\begin{array}{cc}a&0\\0&0\end{array}\right ])1_{\mathcal A}-1_{\mathcal B}\delta(\left
[ \begin{array}{cc}0&0\\0&b\end{array}\right ])1_{\mathcal B}\\
&= (\delta(\left [ \begin{array}{cc}a&0\\0&0\end{array}\right
])-1_{\mathcal A}\delta(\left [
\begin{array}{cc}a&0\\0&0\end{array}\right
])1_{\mathcal A})+\delta(\left [ \begin{array}{cc}0&m\\0&0\end{array}\right ])\\
&+ (\delta(\left [ \begin{array}{cc}0&0\\0&b\end{array}\right
])-1_{\mathcal B}\delta(\left [
\begin{array}{cc}0&0\\0&b\end{array}\right
])1_{\mathcal B})\\
&=1_{\mathcal B}\delta(\left [
\begin{array}{cc}1_{\mathcal A}&0\\0&0
 \end{array}\right ])1_{\mathcal A}\sigma(\left[\begin{array}{cc}a&0\\0&0
 \end{array}\right])+ 1_{\mathcal B}\delta(\left [ \begin{array}{cc}1_{\mathcal A}&0\\
 0&0\end{array}\right ])1_{\mathcal A}
\sigma(\left[\begin{array}{cc}0&m\\0&b\end{array}\right])\\
&-\tau(\left
[\begin{array}{cc}a&m\\0&0\end{array}\right])1_{\mathcal
B}\delta(\left [
\begin{array}{cc}1_{\mathcal A}&0\\0&0\end{array}\right ])1_{\mathcal A} -\tau(\left[
\begin{array}{cc}0&0\\0&b\end{array}\right])
 1_{\mathcal B}\delta(\left [ \begin{array}{cc}1_{\mathcal A}&0\\0&0\end{array}\right
  ])1_{\mathcal A} \\
&= -\delta_{1_{\mathcal B}\delta(\left [
\begin{array}{cc}1_{\mathcal A}&0\\0&0\end{array}\right ])1_{\mathcal A}}(\left [
\begin{array}{cc}
a&m\\0&b\end{array}\right ]).
\endaligned\]
We therefore have $\delta-D\in N^1_{(\sigma, \tau)}({\mathcal
T},\mathcal X)$, and so $\delta\in N^1_{(\sigma, \tau)}({\mathcal
T},\mathcal X).$

Conversely, let $\delta\in N^1_{(\sigma , \tau)}({\mathcal
T},\mathcal X)$. Then there exists $x\in\mathcal  X$ such that
$$\delta(\left [
\begin{array}{cc}a&m
\\0&b\end{array}\right ])= \tau(\left [ \begin{array}{cc}a&m\\0&b\end{array}
\right ])x-x\sigma(\left [ \begin{array}{cc}a&m\\0&b
\end{array}\right ]).$$
Hence
\begin{eqnarray*}
\delta_{\mathcal A}(a)&=&1_{\mathcal A}\delta(\left [
\begin{array}{cc}a&0\\0&0\end{array}\right ])1_{\mathcal A}\\
&=&1_{\mathcal A}(\tau(\left [
\begin{array}{cc}a&0\\0&0\end{array}\right ])x-x\sigma(\left [ \begin{array}{cc}
a&0\\0&0\end{array}\right ]))1_{\mathcal A} \\
&=&\tau(\left [
\begin{array}{cc}a&0\\0&0\end{array}\right ])1_{\mathcal A}x1_{\mathcal A}-1_{\mathcal A}x1_{\mathcal A}\sigma(\left [
\begin{array}{cc}a&0\\0&0\end{array}\right ])\\
&=&\delta_{1_{\mathcal A}x1_{\mathcal A}}(a).
\end{eqnarray*}
Similarly, $\delta_{\mathcal B}(b)=\delta_{1_{\mathcal
B}x1_{\mathcal B}}(b)$. Hence $\delta_{\mathcal A}$ and
$\delta_{\mathcal B}$ are inner and so $\delta\in {\rm ker}\rho$.\\
Thus $N^1_{(\sigma-\tau)}({\mathcal T},\mathcal X)={\rm ker}\rho$.\\
We conclude that
$$H^1_{(\sigma, \tau)}({\mathcal T},\mathcal X)=\frac{Z^1_{(\sigma, \tau)}
({\mathcal T},\mathcal X)}{N^1_{(\sigma, \tau)}({\mathcal T},
\mathcal X)}=\frac{Z^1_{(\sigma, \tau)}({\mathcal T},\mathcal
X)}{{\rm ker}\rho} =H^1_{(\sigma, \tau)}(\mathcal A,\mathcal
X_{\mathcal A\mathcal A})\oplus H^1_{(\sigma, \tau)}(\mathcal
B,\mathcal X_{\mathcal B\mathcal B}).$$
\end{proof}

\begin{corollary} $H^1_{(\sigma, \tau)}({\rm Tri}(\mathcal A,\mathcal M, \mathcal B),\mathcal M)=0$.
\end{corollary}
\begin{proof} With $\mathcal X=\mathcal M$ we have
$$H^1_{(\sigma, \tau)}({\rm Tri}(\mathcal A, \mathcal M, \mathcal B),\mathcal M)
=H^1_{(\sigma, \tau)}(\mathcal A,0)\oplus H^1_{(\sigma,
\tau)}(\mathcal B,0)=0.$$
\end{proof}

\begin{example} $H^1_{(\sigma, \tau)}({\rm Tri}(\mathcal A, \mathcal A,\mathcal A),\mathcal A)=0$.
\end{example}

\begin{example} Let $\mathcal L$ be a left Banach $\mathcal A$-module. Then
$H^1_{(\sigma, \tau)}({\rm Tri}(\mathcal  A,\mathcal L,{\mathbb
C}),\mathcal L)=0$.
\end{example}

\begin{corollary}
$H^1_{(\sigma-\tau)}({\rm Tri}(\mathcal A,  \mathcal
M, \mathcal B),\mathcal A)=H^1_{(\sigma, \tau)}(\mathcal A, \mathcal A)$.
\end{corollary}

\begin{proof} With $\mathcal X=\mathcal A$, we have $\mathcal X_{\mathcal A\mathcal B}=0,
\mathcal X_{\mathcal A\mathcal A}=\mathcal A$ and $\mathcal
X_{\mathcal B\mathcal B}=0$. It then follows from Theorem
\ref{main}, $H^1_{(\sigma, \tau)}({\rm Tri}(\mathcal A, \mathcal
M, \mathcal B),\mathcal A)=H^1_{(\sigma, \tau)}(\mathcal
A,\mathcal A)\oplus H^1_{(\sigma, \tau)}(\mathcal
B,0)=H^1_{(\sigma, \tau)}(\mathcal A, \mathcal A)$.
\end{proof}

\begin{example} If $\mathcal A$ is a hyperfinite von Neumann algebra and
$\mathcal B$ is an arbitrary unital Banach module, then
$H^1_{(\sigma, \tau)}({\rm Tri}(\mathcal A, \mathcal M, \mathcal
B), \mathcal A)=H^1_{(\sigma, \tau)}(\mathcal A,\mathcal A)=0$, if
$\sigma$ and $\tau$ are ultra-weak automorphisms (see Corollary
3.4.6 of \cite{S-S}).
\end{example}

\section{$(\sigma,\tau)$-weak amenability of triangular Banach algebras}

With simple calculation we can observe that if $\mathcal
X=\mathcal T^*$  considered as $\mathcal T$-bimodule, then
$\mathcal X_{\mathcal A\mathcal A}=\mathcal A^*, \mathcal
X_{\mathcal B\mathcal B}=\mathcal B^*$ and $\mathcal X_{\mathcal
A\mathcal B}=0$. Therefore by Theorem \ref{main} we can conclude
the following

\begin{theorem}
Let $\mathcal A,\mathcal B$ be unital Banach algebras  and
$\mathcal M$ be a unital Banach $\mathcal A-\mathcal B$-bimodule.
Then
$$H^1_{(\sigma,\tau)}({\mathcal
T},\mathcal T^*)=H^1_{(\sigma,\tau)}(\mathcal A,\mathcal
A^*)\oplus H^1_{(\sigma,\tau)}(\mathcal B,\mathcal B^*).$$
\end{theorem}

\begin{corollary}
Let $\mathcal A,\mathcal B$ be unital Banach algebras  and
$\mathcal M$ be an unital Banach $\mathcal A-\mathcal B$-bimodule.
The triangular Banach algebra $\mathcal T=Tri(\mathcal A,
\mathcal M, \mathcal B)$ is $(\sigma,\tau)$-weak amenable if and
only if $\mathcal A$ and $\mathcal B$ are both
$(\sigma,\tau)$-weak amenable.
\end{corollary}

By induction one can easily prove the following proposition

\begin{lemma}\label{nterm}
Suppose that $\mathcal A,\mathcal B$ are unital Banach algebras
and $\mathcal M$ is a unital Banach $\mathcal A-\mathcal
B$-bimodule. If $\mathcal X=\mathcal T^{(2n)}$ then
$$\mathcal
X_{\mathcal A\mathcal A}=\mathcal A^{(2n)}, \mathcal X_{\mathcal
B\mathcal B}=\mathcal B^{(2n)}, \mathcal X_{\mathcal A\mathcal
B}=\mathcal M^{(2n)}, \mathcal X_{\mathcal B\mathcal A}=0.$$ Also
if $\mathcal X=\mathcal T^{(2n-1)}$ then $$\mathcal X_{\mathcal
A\mathcal A}=\mathcal A^{(2n-1)}, \mathcal X_{\mathcal B\mathcal
B}=\mathcal B^{(2n-1)}, \mathcal X_{\mathcal A\mathcal B}=0,
\mathcal X_{\mathcal B\mathcal A}=\mathcal M^{(2n-1)}.$$
\end{lemma}

Now by Lemma \ref{nterm} and Theorem \ref{main}, we immediately
 obtain the next result.

\begin{proposition}
Let $\mathcal A,\mathcal B$ be unital Banach algebras and
$\mathcal M$ be a  unital Banach $\mathcal A-\mathcal B$-bimodule.
Then for all positive integers $n\in\mathbb{N}$,
$$H^1_{(\sigma,\tau)}({\mathcal
T},\mathcal T^{(2n-1)})=H^1_{(\sigma,\tau)}(\mathcal A,\mathcal
A^{(2n-1)})\oplus H^1_{(\sigma,\tau)}(\mathcal B,\mathcal
B^{(2n-1)}).
$$
\end{proposition}

\section{$(\sigma,\tau)$-amenability of triangular Banach algebras}

In this section, by using some ideas of \cite{M-S-Y} we investigate
$(\sigma,\tau)$-amenability of triangular Banach algebra $\mathcal
T=Tri(\mathcal A, \mathcal M, \mathcal B)$. We shall assume that the
homomorphisms $\sigma, \tau$ on $\mathcal T$ have properties
asserted  in (\ref{equ1}) and (\ref{equ2}). We need some general
observation concerning $(\sigma, \tau)$-amenability of Banach
algebras. The first is an easy consequence of the definition of
$(\sigma, \tau)$-amenability.
\begin{proposition}\cite[Proposition 3.3]{M-N}\label{old1}
Let ${\mathcal A},{\mathcal   B}$ be Banach algebras and
$\sigma,\sigma'$ be continuous endomorphisms of ${\mathcal    A}$
and $\tau, \tau'$ be continuous homomorphisms of ${\mathcal B}$. If
there is a continuous homomorphism $\varphi:{\mathcal
A}\longrightarrow {\mathcal B}$ such that $\varphi({\mathcal A})$ is
a dense subalgebra of ${\mathcal B}$ and $\tau\varphi=\varphi\sigma$
and $\tau'\varphi=\varphi\sigma'$, then $(\sigma,
\sigma')$-amenability of ${\mathcal   A}$ implies $(\tau, \tau')
$-amenability of ${\mathcal B}$.
\end{proposition}
Now, suppose that $\mathcal A$ is a Banach algebra,
$\tau,\sigma:{\mathcal A}\longrightarrow {\mathcal A}$ are two
continuous endomorphisms, and ${\mathcal   I}$ is a closed ideal of
${\mathcal   A}$ such that $\sigma({\mathcal I})\subseteq {\mathcal
I}, \tau({\mathcal I})\subseteq {\mathcal   I}$. Then the map
$\widehat\tau,\widehat\sigma:\frac{\mathcal   A}{\mathcal
I}\longrightarrow \frac{\mathcal   A}{\mathcal   I}$ can be defined
by $\widehat\sigma(a+{\mathcal   I})=\sigma(a)+{\mathcal I},
\widehat\tau(a+{\mathcal   I}) =\tau(a)+{\mathcal   I}$. It is not
hard to show the following propositions.
\begin{proposition}\cite[Proposition 3.1]{M-N}\label{old2}
Let ${\mathcal   I}, \sigma, \tau$ be as above. If ${\mathcal A}$ is
$(\sigma,\tau)$-amenable then $\frac{\mathcal   A}{\mathcal I}$ is
$(\widehat\sigma,\widehat\tau)$-amenable.
\end{proposition}
\begin{proposition}\cite[Proposition 3.2]{M-N}\label{old3}
Let ${\mathcal   I}, \sigma,\tau$ be as above and let $\sigma, \tau$
be idempotent homomorphisms. If ${\mathcal   I}$ is
$(\sigma,\tau)$-amenable and $\frac{{\mathcal   A}}{{\mathcal I}}$
is $(\widehat\sigma,\widehat\tau)$-amenable, then ${\mathcal A}$ is
$(\sigma,\tau)$-amenable.
\end{proposition}

We now extend Theorem 4.1 of \cite{M-S-Y} as follows

\begin{proposition}
Let $\sigma$ and $\tau$ be two continuous idempotent homomorphisms
on triangular Banach algebra $\mathcal T=Tri(\mathcal A, 0, \mathcal
B)$. The triangular Banach algebra $\mathcal T$ is
$(\sigma,\tau)$-amenable if and only if $\mathcal A$ and $\mathcal
B$ are $(\sigma,\tau)$-amenable.
\end{proposition}
\begin{proof}
At first suppose that $\mathcal A, \mathcal B$ are
$(\sigma,\tau)$-amenable. It is easy to see that $ \left [
\begin{array}{cc}\mathcal A&0\\0&0 \end{array}\right ]$ is a
closed ideal of $ \left [
\begin{array}{cc}\mathcal A&0\\0&\mathcal B \end{array}\right
]$ and $ \left [
\begin{array}{cc}\mathcal A&0\\0&\mathcal B \end{array}\right ]/\left [
\begin{array}{cc}\mathcal A&0\\0&0 \end{array}\right ]\simeq \left [
\begin{array}{cc}0&0\\0&\mathcal B \end{array}\right ]$. Since $\mathcal A$ is
$(\sigma,\tau)$-amenable therefore $ \left [
\begin{array}{cc}\mathcal A&0\\0&0 \end{array}\right ]$ is $(\sigma,\tau)$-amenable.
Let $\varphi: \left [
\begin{array}{cc}0&0\\0&\mathcal B \end{array}\right ] \to \left [
\begin{array}{cc}\mathcal A&0\\0&\mathcal B \end{array}\right ]/ \left [
\begin{array}{cc}\mathcal A&0\\0&0 \end{array}\right ] $  be the natural isomorphism. Then
$\varphi\tau=\widehat{\tau}\varphi$ and
$\varphi\sigma=\widehat{\sigma}\varphi$. By Proposition
\ref{old1} $(\sigma,\tau)$-amenability of $\left [
\begin{array}{cc}0&0\\0&\mathcal B \end{array}\right ]$ implies the
$(\widehat{\sigma},\widehat{\tau})$-amenability of $\left [
\begin{array}{cc}\mathcal A&0\\0&\mathcal B \end{array}\right ]/ \left [
\begin{array}{cc}\mathcal A&0\\0&0 \end{array}\right ]$. Thus by utilizing Proposition \ref{old3}, we deduce
the  $(\sigma,\tau)$-amenability of the Banach algebra $\mathcal T$.

For the converse, suppose that $\mathcal T$ is
$(\sigma,\tau)$-amenable. It is obvious that $\left [
\begin{array}{cc}\mathcal A&0\\0&0 \end{array}\right
]$ is a closed ideal of $\mathcal T$. By Proposition \ref{old2},
$\left [
\begin{array}{cc}\mathcal A&0\\0&\mathcal B \end{array}\right
]/\left [
\begin{array}{cc}\mathcal A&0\\0&0 \end{array}\right
]$ is $(\widehat{\sigma},\widehat{\tau})$-amenable. One can easily
observe that there exists the natural isomorphism $\varphi:\left [
\begin{array}{cc}\mathcal A&0\\0&\mathcal B \end{array}\right
]/\left [
\begin{array}{cc}\mathcal A&0\\0&0 \end{array}\right] \to \left [
\begin{array}{cc}0&0\\0&\mathcal B \end{array}\right ]$ and that  $\varphi\widehat{\sigma}=\sigma\varphi$ and
$\varphi\widehat{\tau}=\tau\varphi$. Therefore, by Proposition
\ref{old1}, $\left [
\begin{array}{cc}0&0\\0&\mathcal B \end{array}\right ]$, that is $\mathcal B$ ,is $(\sigma,\tau)$-amenable.
Similarly one can prove the $(\sigma,\tau)$-amenability of $\mathcal
A$.
\end{proof}
\begin{theorem}
Let $\sigma$ and $\tau$ be two continuous idempotent homomorphisms
on triangular Banach algebra $\mathcal T=Tri(\mathcal A,\mathcal M,
\mathcal B)$.  If the triangular Banach algebra $\mathcal T$ is
$(\sigma,\tau)$-amenable then $\mathcal A$ and $\mathcal B$ are
$(\sigma,\tau)$- amenable. In particular, $\sigma$-amenability of
$\mathcal T$  implies $\sigma(\mathcal M)=\{0\}$
\end{theorem}
\begin{proof}
Suppose that $\mathcal T=\left [
\begin{array}{cc}\mathcal A&\mathcal M\\0&\mathcal B \end{array}\right
]$ is $(\sigma,\tau)$-amenable. Clearly, $\left [
\begin{array}{cc}\mathcal A&\mathcal M\\0&0 \end{array}\right
]$ is a closed ideal of $\mathcal T$. Therefore, by Proposition
\ref{old2}, $\left [
\begin{array}{cc}\mathcal A&\mathcal M\\0&\mathcal B \end{array}\right
]/\left [\begin{array}{cc}\mathcal A&\mathcal
M\\0&0\end{array}\right ]$ is
$(\widehat{\sigma},\widehat{\tau})$-amenable. Also there exists
the natural isomorphism $\varphi:\left [
\begin{array}{cc}\mathcal A&\mathcal M\\0&\mathcal B \end{array}\right
]/\left [\begin{array}{cc}\mathcal A&\mathcal
M\\0&0\end{array}\right ] \to \left [
\begin{array}{cc}0&0\\0&\mathcal B \end{array}\right
]$ such that $\varphi\widehat{\sigma}=\sigma\varphi$ and
$\varphi\widehat{\tau}=\tau\varphi$. Hence $\left [
\begin{array}{cc}0&0\\0&\mathcal B \end{array}\right
]$ is $(\sigma,\tau)$-amenable. Similarly one can prove the
$(\sigma,\tau)$-amenability of $\mathcal A$.

Now suppose that the triangular Banach algebra $\mathcal T$ is
$\sigma$-amenable. Set $\mathcal X=\left [
\begin{array}{cc}\mathcal A^*&\mathcal M^*\\0&\mathcal B^* \end{array}\right
]$. The vector space $\mathcal X$ is a Banach space under the norm
$\|\left [
\begin{array}{cc}f&h\\0&g
\end{array}\right ]\|=\|f\|_{\mathcal A^*}+\|h\|_{\mathcal M^*}+\|g\|_{\mathcal
B^*}$. The space $\mathcal X$ can be regarded as  a Banach $\mathcal
T$-bimodule under the following $\mathcal T$-module actions
\begin{eqnarray*}
\left [
\begin{array}{cc}a&m\\0&b \end{array}\right ]\cdot\left [
\begin{array}{cc}f&h\\0&g \end{array}\right ]&=&\left [
\begin{array}{cc}0&bh\\0&0 \end{array}\right ]\\
\left [
\begin{array}{cc}f&h\\0&g \end{array}\right ]\cdot\left [
\begin{array}{cc}a&m\\0&b \end{array}\right ]&=&\left [
\begin{array}{cc}0&ha\\0&0 \end{array}\right ]
\end{eqnarray*}
where $ \left [
\begin{array}{cc}a&m\\0&b \end{array}\right ]\in\mathcal T$ and $ \left [
\begin{array}{cc}f&h\\0&g \end{array}\right ]\in \mathcal X$.
Therefore   $\mathcal X^*= \left [
\begin{array}{cc}\mathcal A^{**}&M^{**}\\0&B^{**} \end{array}\right
]$  is a dual Banach $\mathcal T$-bimodule.

Let $D:\mathcal T\longrightarrow\left [
\begin{array}{cc}\mathcal A^{**}&M^{**}\\0&B^{**} \end{array}\right
]$ be defined by $D(\left [
\begin{array}{cc}a&m\\0&b \end{array}\right ])=\left [
\begin{array}{cc}0&\widehat{\sigma(m)}\\0&0\end{array}\right ]$.
Now we have
\begin{eqnarray*}
D\Big(\left [
\begin{array}{cc}a_1&m_1\\0&b_1 \end{array}\right ]\left [
\begin{array}{cc}a_2&m_2\\0&b_2 \end{array}\right
]\Big)&=&D\Big(\left [
\begin{array}{cc}a_1a_2&a_1m_2+m_1b_2\\0&b_1b_2 \end{array}\right
]\Big)\\
&=&\left [
\begin{array}{cc}0&\widehat{\sigma(a_1m_2+m_1b_2)}\\0&0 \end{array}\right
]\\
&=&\left [
\begin{array}{cc}0&\sigma(a_1)\widehat{\sigma(m_2)}+\widehat{\sigma(m_1)}\sigma(b_2)\\0&0
\end{array}\right]\\
&=&\left [
\begin{array}{cc}0&\sigma(a_1)\widehat{\sigma(m_2)}\\0&0
\end{array}\right]+\left [
\begin{array}{cc}0&\widehat{\sigma(m_1)}\sigma(b_2)\\0&0
\end{array}\right]\\
&=&\left [
\begin{array}{cc}0&\widehat{\sigma(m_1)}\\0&0
\end{array}\right]\left [
\begin{array}{cc}\sigma(a_2)&\sigma(m_2)\\0&\sigma(b_2) \end{array}\right
]\\
&+& \left [
\begin{array}{cc}\sigma(a_1)&\sigma(m_1)\\0&\sigma(b_1) \end{array}\right
]\left [
\begin{array}{cc}0&\widehat{\sigma(m_2)}\\0&0
\end{array}\right]\\
&=&D\Big(\left [
\begin{array}{cc}a_1&m_1\\0&b_1 \end{array}\right
]\Big)\sigma\Big(\left [
\begin{array}{cc}a_2&m_2\\0&b_2 \end{array}\right ]\Big)\\
&+&\sigma\Big(\left [
\begin{array}{cc}a_1&m_1\\0&b_1 \end{array}\right
]\Big)D\Big(\left [
\begin{array}{cc}a_2&m_2\\0&b_2 \end{array}\right ]\Big).
\end{eqnarray*}
Therefore $D$ is a $\sigma$-derivation. Hence there exists $\left [
\begin{array}{cc}F&H\\0&G \end{array}\right ]\in\mathcal X^*$ such that
\begin{eqnarray*}
\left [
\begin{array}{cc}0&\widehat{\sigma(m)}\\0&0\end{array}\right ]&=&D\Big(\left [
\begin{array}{cc}a&m\\0&b \end{array}\right ]\Big)\\
&=&\sigma\Big(\left [
\begin{array}{cc}a&m\\0&b \end{array}\right ]\Big)\left [
\begin{array}{cc}F&H\\0&G \end{array}\right ]-\left [
\begin{array}{cc}F&H\\0&G \end{array}\right ]\sigma\Big(\left [
\begin{array}{cc}a&m\\0&b \end{array}\right ]\Big)\\
&=&\left [
\begin{array}{cc}0&\sigma(a)H\\0&0 \end{array}\right ]-\left [
\begin{array}{cc}0&H\sigma(b)\\0&0 \end{array}\right ].
\end{eqnarray*}
Thus
\begin{eqnarray}\label{hat}
\widehat{\sigma(m)}=\sigma(a)H-H\sigma(b)
\end{eqnarray}
for all $m\in\mathcal M, a\in\mathcal A$ and $b\in\mathcal B$.
Choosing $a=0, b=0$  in (\ref{hat}) we conclude that
$\widehat{\sigma(m)}=0$ for all $m\in\mathcal M$. Thus
$\sigma(\mathcal M)=\{0\}$.
\end{proof}
\begin{corollary}
Let $\mathcal A, \mathcal B$ be two unital Banach algebra and
$\sigma:\mathcal A\longrightarrow\mathcal A, \tau:\mathcal
B\longrightarrow\mathcal B$ be two continuous  idempotent
homomorphisms. The Banach algebra $\mathcal A\oplus\mathcal B$ is
$\sigma\oplus\tau$-amenable if and only if $\mathcal A$ is
$\sigma$-amenable and $\mathcal B$ is $\tau$-amenable.
\end{corollary}
\begin{proof}
It is easy to see that $\left [
\begin{array}{cc}\mathcal A&0\\0&\mathcal B \end{array}\right
]\simeq\mathcal A\oplus\mathcal B$. Define $\varphi:\left [
\begin{array}{cc}\mathcal A&0\\0&\mathcal B \end{array}\right
]\longrightarrow\left [
\begin{array}{cc}\mathcal A&0\\0&\mathcal B \end{array}\right
]$ via $\varphi\Big(\left [
\begin{array}{cc} a&0\\0& b \end{array}\right
]\Big)=\left [
\begin{array}{cc}\sigma(a)&0\\0&\tau(b) \end{array}\right
]$. Therefore $\left [
\begin{array}{cc}\mathcal A&0\\0&\mathcal B \end{array}\right
]$ is $\varphi$-amenable if and only if $\mathcal A$ and $\mathcal
B$ are both $\varphi$-amenable, and this holds if and only if
$\mathcal A$ is $\sigma$-amenable and $\mathcal B$ is
$\tau$-amenable.
\end{proof}

\end{document}